\documentclass[12pt]{article}

\usepackage[margin=1in]{geometry}
\usepackage{amsmath}
\usepackage{amsfonts}
\usepackage{amsthm}
\usepackage{mathtools}
\usepackage{hyperref}
\usepackage{amssymb}

\newtheorem{thm}{Theorem}[section]
\newtheorem{lem}[thm]{Lemma}
\newtheorem{prop}[thm]{Proposition}
\newtheorem{cor}[thm]{Corollary}

\theoremstyle{definition}
\newtheorem{defn}[thm]{Definition}
\theoremstyle{plain}

\newcommand{\HH}{\mathcal H}
\newcommand{\RR}{\mathbb R}
\newcommand{\ddc}{\sqrt{-1}\partial\bar\partial}

\newcommand{\Pol}{\operatorname{Pol}}
\numberwithin{equation}{section}

\begin{document}

\title{Coerciveness for complex Hessian type equations}
\author{Gao Chen, Junbang Liu, and Sijie Nie}
\date{}
\maketitle

\begin{abstract}
We study complex Hessian type equations associated with right-Noetherian polynomials on compact K\"ahler manifolds of complex dimension $n$. For two such polynomials of the same degree $2\le d\le n$ in strict proper position, we prove that under the integrability condition and the assumption that the admissible space $\HH_2$ is non-empty, smooth admissible solvability for $f_1$, the $d_1$-coerciveness of the associated variational functional $J_1$ on normalized potentials, and the numerical criterion for $f_1$ are equivalent to each other. This extends the existence-coerciveness correspondence of Chen-Ghosh-Wang to the lower-degree setting under weaker polynomial conditions.
\end{abstract}

\tableofcontents

\section{Introduction}
\label{sec:introduction}

This paper continues the work of Chen-Ghosh-Wang \cite{ChenGhoshWangCoerciveness} on the relation between coerciveness and solvability of generalized Monge-Amp\`ere equations. We consider the complex Hessian type equations associated to a right-Noetherian polynomial of degree $d$ on a compact K\"ahler manifold of dimension $n$, with the possibility that $d<n$. There are three related problems: the solvability of the equation, the numerical criterion, and the coerciveness criterion. We show that they are all equivalent to each other.

The motivation comes from the properness problem in K\"ahler geometry. 
A central problem in K\"ahler geometry is to search for a canonical metric in a K\"ahler class. The most important one is a constant scalar curvature K\"ahler (cscK) metric. It is a critical point of the Mabuchi K-energy \cite{MabuchiKEnergy}. The equivalence of existence of the cscK metrics and the coerciveness of the Mabuchi K-energy in $d_1$-distance modulo holomorphic automorphisms was established by Berman-Darvas-Lu \cite{BermanDarvasLuRegularity} and Chen-Cheng \cite{ChenChengExistence} in two directions. By writing the Mabuchi K-energy as the sum of entropy and the $J$-functional, Chen introduced the $J$-functional \cite{ChenLowerBoundMabuchi}, whose critical point equation is the $J$-equation introduced by Donaldson from a moment map picture \cite{DonaldsonMomentMaps}. For the $J$-equation, the equivalence between solvability and coerciveness was proved by Collins-Sz\'ekelyhidi \cite{CollinsSzekelyhidiJFlow}. A further example is the Leung-Yau-Zaslow equation \cite{LeungYauZaslowFourierMukai}, also called the deformed Hermitian-Yang-Mills equation: the existence-coerciveness equivalence was proved by Chu-Lee \cite{ChuLeeHypercriticalDHYM} in the hypercritical phase when the complex dimension is at least $4$. These results suggest that the variational characterization should be studied for broader classes of complex Hessian type equations. 

Chen-Ghosh-Wang developed such a theory for degree-$n$ right-Noetherian polynomials, extending the coerciveness results to the supercritical LYZ equation. The setting of such polynomials was introduced by Lin \cite{LinInverseSigmaConvexity,LinInverseSigmaSolvability}, and has been studied in more detail by Fang-Ma \cite{FangMaGarding, FangMaIdealGarding, FangMaSpectral}.   In the theory of Chen-Ghosh-Wang, they use two polynomials $f_1,f_2$ in the strongly strict proper position, and with the ideal proper-position $f_2\prec \prec f_1$. The ideal proper position implies the concavity of $g_1/g_2$, which is used to establish the higher-order regularity of the flow $\dot \varphi=\frac{G_1(\chi_\varphi)}{G_2(\chi_\varphi)}$. In the present paper, we use the numerical criterion as a middle step to prove that coerciveness implies solvability, and we only need the strict proper-position condition.

We describe the setting now. Let $(X,\omega)$ be a compact connected K\"ahler manifold of complex dimension $n$, and let $\chi$ be a smooth real closed $(1,1)$-form. Write $\chi_\varphi:=\chi + \ddc \varphi$. For $i=1,2$, we define the polarizations $g_i$, $G_i$ of $f_i$ as in \eqref{eq:polynomials-forms}. Following the variational setup in \cite{ChenGhoshNieGeodesics, ChenGhoshWangCoerciveness, GhoshWangXuGeodesics}, we define the admissible space $\HH_2$ by \eqref{eq:H2-definition} with the assumption that $0\in \mathcal H_2$, and define the $d_p$-distance by \eqref{eq:d1-definition}. We define the functionals $I_i$, for $i=1,2$, by the explicit formula
\eqref{eq:energies} and the variational formula \eqref{eq:energy-variation}.

\begin{thm}
\label{thm:main}
Assume that $f_1,f_2$ are monic right-Noetherian polynomials with the same degree $2\le d\le n$, and $f_2\prec f_1$ strictly, with $0\in\HH_2$. Assume that $V_1=\int_X G_1(\chi)=0$.
Then the followings are equivalent:
\begin{enumerate}
\item There is a smooth $\varphi_*\in\HH_2$ such that
      $g_1(\lambda(\omega^{-1}\chi_{\varphi_*}))=0$.
\item There are constants $\delta>0$ and $C\ge0$ such that
\begin{equation}\label{eq:main-d1-coercivity}
 J_1(\varphi)\ge\delta d_1(0,\varphi)-C
 \quad\text{for all }\varphi\in\HH_2\text{ with }I_2(\varphi)=0.
\end{equation}
\item For every irreducible subvariety $Y$ with codimension $1 \le r\le d$, 
\begin{equation}
 \int_Y G_1^{(r)}(\chi)>0.
\end{equation}
\end{enumerate}
\end{thm}

 The main difference in the lower-degree case compared to \cite{ChenGhoshWangCoerciveness} is about the direction in which coerciveness implies solvability.  For any subvariety $Y$, we construct a ray in the admissible space $\mathcal H_2$. The computation of the asymptotic slope of $I_2$ implies the numerical criterion; see Section~\ref{sec:numerical}. A similar calculation was carried out in the special case of the $J$-equation and when $Y$ is a divisor in \cite{LiuSemistableJEquation}. To apply the numerical criterion, we establish a connected component principle for the solvability, based on Fu-Li-Zhang-Zhang's \cite{FuLiZhangZhangNumerical} extension of Chen-Nie-Xu \cite{ChenNieXuNumericalCriterion} to the K\"ahler setting, and Fang-Ma's subsolution to solution theorem \cite{FangMaSpectral}. This completes the equivalence without using the parabolic flow, and thus we can avoid the higher order regularity issue and remove the condition of ideal proper position. 
 
 The paper is organized as follows. Section~\ref{sec:preliminaries} discusses the Fang-Ma-G\aa rding framework and recalls the variational setup of \cite{ChenGhoshNieGeodesics, ChenGhoshWangCoerciveness, GhoshWangXuGeodesics}. 
The direction $(1) \implies (2)$ is proved in Section~\ref{sec:solvability-coerciveness}. The direction $(2)\implies (3)$ is proved in Section~\ref{sec:numerical}. The direction $(3)\implies (1)$ is proved in Section~\ref{sec:continuity}.

\paragraph{\textbf{Declaration on the use of AI.}}
ChatGPT 6 Astra was used to help implement the authors' ideas, and we take responsibility for all mathematical arguments and conclusions. We first asked AI to use \cite{ChenGhoshWangCoerciveness} to give a proof. This attempt did not yield a complete proof because of the difficulty of finding an admissible smooth approximation to currents. Then we asked AI to use \cite{FuLiZhangZhangNumerical}. This attempt also did not yield a complete proof because of the difficulty of computing the mixed integrals directly from coerciveness. Finally, inspired by \cite{ChenChengExistence}, we came up with the idea of using the connected component argument, and AI produced the first draft of this paper. After that, we identified arguments already present in the literature and added the corresponding citations. We also added missing details and simplified the notation. Finally, AI was used to assist with checking for typos and correcting grammar.

\section{Preliminaries}
\label{sec:preliminaries}
We recall the polynomial conventions and geometric setup used in
\cite{FangMaGarding, FangMaIdealGarding, ChenGhoshNieGeodesics, ChenGhoshWangCoerciveness, GhoshWangXuGeodesics}.
Throughout, $(X,\omega)$ is a compact connected K\"ahler manifold of complex dimension $n$, and $\chi$ is a smooth real closed $(1,1)$-form.

We first recall the definition of right-Noetherian polynomials:

\begin{defn}[Definition~2.1 of \cite{LinInverseSigmaConvexity}]
\label{def:right-Noetherian}
Let $f\in\RR[x]$ be a polynomial with degree $1\le d\le n$ and positive leading coefficient. We call $f$
\emph{right-Noetherian} if the largest real root $r(f^{(k)})$ of $f^{(k)}$ exists for every $0\le k\le d-1$ and
\begin{equation}\label{eq:monotone-root-sequence}
 r(f)\ge r(f')\ge\cdots\ge r(f^{(d-1)}).
\end{equation}
For $d\ge2$, $f$ is \emph{strictly right-Noetherian} if the first inequality is strict, and \emph{strongly strictly right-Noetherian} if every inequality is strict.
\end{defn}
Remark that the multiplication by a positive constant doesn't change the roots, and only changes the functional and distance by a multiple, so we are free to use inconsistent notations as long as they can be fixed by such a multiplication.

Next we recall the definition of proper position in \cite{FangMaGarding}.

\begin{defn}[Theorem~10.12 of \cite{FangMaGarding}]
\label{def:proper-position}
Let $f_1,f_2$ be monic right-Noetherian polynomials of degree $d$. We say that $f_2$ and $f_1$ are \emph{in proper position}, and write $f_2\prec f_1$, if, for every $0\le k\le d-1$,
\begin{equation}\label{eq:univariate-proper-position}
 r(f_2^{(k)})\le r(f_1^{(k)}),\qquad \text{and }
 f_1^{(k)}(x)\le 0\quad\text{for }r(f_2^{(k)})\le x<r(f_1^{(k)}).
\end{equation}
\end{defn}

We also need a strict version as in \cite{ChenGhoshWangCoerciveness}:

\begin{defn}[Section~3 of \cite{ChenGhoshWangCoerciveness}]
\label{def:strongly-strict-proper-position}
Assume $f_2\prec f_1$. We say that this
proper position is \emph{strict} if
\begin{equation}\label{eq:strict-proper-position-pole}
 \lim_{x\downarrow r(f_2)}\frac{f_1(x)}{f_2(x)}=-\infty.
\end{equation}
We say that $f_2\prec f_1$ \emph{strongly strictly} if
$f_2^{(k)}\prec f_1^{(k)}$ strictly for every $0\le k\le d-1$.
\end{defn}

The polynomial assumptions used below are
\begin{equation}\label{eq:polynomial-assumptions}
 \begin{gathered}
 f_1,f_2\text{ are monic right-Noetherian polynomials},\\
 \deg f_1=\deg f_2=d,\qquad 2\le d\le n,\qquad
 f_2\prec f_1\text{ strictly}.
 \end{gathered}
\end{equation}
We have removed the \emph{ideal proper position} condition defined in \cite{FangMaIdealGarding} and used in \cite{ChenGhoshWangCoerciveness}, and we have relaxed the strongly strict proper position condition in \cite{ChenGhoshWangCoerciveness} to just strict proper position.

Next we recall the polarization. For $1\le k\le n$, let
\[
 \sigma_k(\lambda_1,\ldots,\lambda_n)
 =\sum_{1\le j_1<\cdots<j_k\le n}
       \lambda_{j_1}\cdots\lambda_{j_k}
\]
be the $k$th elementary symmetric polynomial in $n$ variables, with
$\sigma_0=1$. For our two polynomials, set
\begin{equation}\label{eq:polynomials-forms}
\begin{aligned}
 f_i(x)&=\sum_{k=0}^d a_k^{(i)}\binom nk x^k,\qquad a_d^{(i)}=\binom nd^{-1},\\
 g_i(\lambda)&=\Pol_n(f_i)(\lambda)
 =\sum_{k=0}^d a_k^{(i)}\sigma_k(\lambda),  \qquad G_i(\chi_\varphi)=\sum_{k=0}^d a_k^{(i)} \binom nk
                 \chi_\varphi^k\wedge\omega^{n-k}.
\end{aligned}
\end{equation}

For $0\le r\le d$, we also define the $r$th formal derivative of $G_i$ with respect to $\chi_\varphi$ by
\[
G_i^{(r)}(\chi_\varphi)
 =\sum_{k=r}^d a_k^{(i)} \frac{k!}{(k-r)!}\binom nk
       \chi_\varphi^{k-r}\wedge\omega^{n-k} =\frac{n!}{(n-r)!}\sum_{k=0}^{d-r} a_{k+r}^{(i)}
       \binom {n-r}{k}\chi_\varphi^{k}\wedge\omega^{n-r-k}.
\]
This is an $(n-r,n-r)$-form. We also write $G_i'(\chi_\varphi)=G_i^{(1)}(\chi_\varphi)$.

The following theorem is one of the main theorems in \cite{FangMaIdealGarding}:
\begin{thm}[Theorem 8.2 of \cite{FangMaIdealGarding}]
If $f$ is monic right-Noetherian, then $g=\Pol_n(f)$ is ideal
Fang-Ma-G\aa rding. In other words, there exists a unique convex connected component $\Upsilon_g$ of $\{g>0\}$ which passes the positive-ray test
\begin{equation}\label{eq:PRT}
 \Upsilon_g+\RR_{\ge0}^n\subseteq\Upsilon_g.
\end{equation}
The same conclusions hold for every nonzero partial derivative
$\partial_I g:=(\prod_{j\in I}\partial_j)g$, with
$\partial_{\varnothing}g=g$. Moreover,
$\Upsilon_{\partial_I g}\subseteq\Upsilon_{\partial_J g}$
whenever $I\subseteq J\subseteq\{1,\ldots,n\}$ and
$\partial_Jg\not\equiv0$. For a positive constant polynomial $c$,
we define $\Upsilon_c=\RR^n$.
\end{thm}
Next, we recall the following key properties of proper position.
\begin{thm}[Theorem 6.14 of \cite{FangMaGarding} and Proposition 3.1 of \cite{ChenGhoshWangCoerciveness}]
$f_2\prec f_1$ implies
\begin{equation}\label{eq:polarized-proper-position}
 \Upsilon_{g_1}=\Upsilon_{g_2}\cap\{g_1>0\},\qquad
 \partial_j\!\left(\frac{g_1}{g_2}\right)
 =\frac{g_2\partial_jg_1-g_1\partial_jg_2}{g_2^2}\ge0
 \quad\text{on }\Upsilon_{g_2}.
\end{equation}
$f_2\prec f_1$ strictly implies $\partial_j(\frac{g_1}{g_2})>0$ on $\Upsilon_{g_2}$.
\end{thm}

Fang-Ma \cite{FangMaGarding} have also proved that the dominance condition $f_2 \triangleleft f_1$ is preserved under polarization. It was also proved by Lin \cite{LinInverseSigmaConvexity} when $d=n$.
\begin{thm}[Theorem 10.10 and Definition 10.8 of \cite{FangMaGarding}]
\label{thm:Dominance}
Suppose that $f_1, f_2$ are monic right-Noetherian polynomials of degrees $d_2\le d_1 \le n$ and $g_1, g_2$ are their polarizations. If for every $0\le k\le d_2-1$,
\[
 r(f_2^{(k)})\le r(f_1^{(k)}),
\]
then
\[
 \Upsilon_{\partial_I g_1}\subseteq\Upsilon_{\partial_I g_2}.
\]
This holds for all $I\subset\{1,...,n\}$ such that $\partial_I g_2$ is not zero. In particular, $I=\emptyset$ means that $\Upsilon_{g_1}\subset\Upsilon_{g_2}$.
\end{thm}

Now we use Proposition 3.1 of \cite{ChenGhoshWangCoerciveness} to prove several key properties of strict proper positions:
\begin{lem}
\label{lem:strict-RN-and-inclusion}
If $f_1, f_2$ are right-Noetherian with degree $2\le d\le n$, and $f_2\prec f_1$ strictly, then $f_1$ is strictly right-Noetherian, and $\bar\Upsilon_{g_1}\subset\Upsilon_{g_2}$.
\end{lem}
\begin{proof}
By Proposition 3.1 of \cite{ChenGhoshWangCoerciveness}, the ratio $\frac{f_1}{f_2}$ tends to $-\infty$ at $r(f_2)$ and to $1$ at infinity. In particular, $r(f_1)>r(f_2)$. At $r(f_1)$, using Remark 3.2 of \cite{ChenGhoshWangCoerciveness}, 
\[f_1'(r(f_1))=f_2(r(f_1))(\frac{f_1}{f_2})'(r(f_1))>0.\]
Thus $f_1$ is strictly right-Noetherian. 

Now we consider the polynomial $f_1+\varepsilon$. All its higher derivatives are the same as $f_1$, and for $\varepsilon$ small, we can make sure that $r(f_1+\varepsilon)>r(f_1')$ and $r(f_1+\varepsilon)>r(f_2)$ still holds. Then $f_1+\varepsilon$ is still right-Noetherian, and $f_2\triangleleft f_1+\varepsilon$. By Theorem~\ref{thm:Dominance},
$\bar\Upsilon_{g_1}\subset\Upsilon_{g_1+\varepsilon}\subset\Upsilon_{g_2}$.
\end{proof}

Following the variational setup in \cite{ChenGhoshNieGeodesics, ChenGhoshWangCoerciveness, GhoshWangXuGeodesics}, we use the admissible space 
\begin{equation}\label{eq:H2-definition}
 \HH_2=\{\varphi\in C^\infty(X,\RR):
       \lambda(\omega^{-1}\chi_\varphi)\in\Upsilon_{g_2}
       \text{ at every point of }X\},
\end{equation}
with the assumption that $0\in \mathcal H_2$. We want to solve the equation $\chi_\varphi\in \partial\Upsilon_{g_1}$ within $\HH_2$. So $\bar\Upsilon_{g_1}\subset\Upsilon_{g_2}$ is a natural condition. The strict right-Noetherian of $f_1$ is also necessary for the analysis in \cite{FangMaSpectral}. That's why we have to assume that $f_2\prec f_1$ strictly.

We define the $d_p$-distance, $1\le p < \infty$, by
\begin{equation}\label{eq:d1-definition}
 d_p(\varphi_0,\varphi_1)
 =\inf_{\varphi_t}\int_0^1 (\int_X
          |\dot\varphi_t|^p G_2(\chi_{\varphi_t}))^{1/p}\,dt,
\end{equation}
where the infimum is over piecewise smooth paths in $\HH_2$
with these endpoints.

Following \cite[Section~7, equation~(10)]{ChenGhoshNieGeodesics},
we define the functionals $I_i$, for $i=1,2$, by the explicit formula
\begin{equation}\label{eq:energies}
 I_i(\varphi)
 =\sum_{k=0}^d a_k^{(i)}\binom nk
       \sum_{j=0}^k\frac{1}{j+1}\binom kj
       \int_X\varphi\,(\ddc\varphi)^j
                    \wedge\chi^{k-j}\wedge\omega^{n-k},
 \qquad J_1=-I_1.
\end{equation}
For a smooth path $\varphi_t\in\HH_2$, we have
\begin{equation}\label{eq:energy-variation}
 \frac{d}{dt}I_i(\varphi_t)
     =\int_X\dot\varphi_t\,G_i(\chi_{\varphi_t}).
\end{equation}

We define the cohomological constants by
\[
 V_i=\int_XG_i(\chi), \qquad v_i=\frac{V_i}{\int_X\omega^n}.
\]
We require that $V_1=0$.

Next, we recall the following theorems about $\rho$-geodesics proved in \cite{GhoshWangXuGeodesics}:

\begin{thm}[Theorem 1.1, Lemma 3.1, and Proposition 6.1 of \cite{GhoshWangXuGeodesics}]
\label{prop:geodesics}
For every $\varphi_0,\varphi_1\in\HH_2$, there are constants
$\rho_0,C>0$ such that, for every $0<\rho\le\rho_0$, there is a
smooth path $\varphi_{\rho,t}\in\HH_2$, $0\le t\le1$, with these
endpoints. Its defect
\begin{equation}\label{eq:geodesic-defect}
 \mathcal D_\rho(t)
 :={}\ddot\varphi_{\rho,t}G_2(\chi_{\varphi_{\rho,t}})-\sqrt{-1}\partial\dot\varphi_{\rho,t}
       \wedge\bar\partial\dot\varphi_{\rho,t}
       \wedge G_2'(\chi_{\varphi_{\rho,t}})
\end{equation}
satisfies
\begin{equation}\label{eq:geodesic-defect-bound}
 0\le\mathcal D_\rho(t)\le C\rho^{1/d}\omega^n.
\end{equation}
\end{thm}

They have also proved the following theorem:

\begin{thm}[Theorem 1.3 of \cite{GhoshWangXuGeodesics}]
\label{ass:distance}
For every $\varphi\in\HH_2$,
\begin{equation}\label{eq:CGN}
 d_1(0,\varphi)\ge
 \max\left\{
       \int_X\varphi_-G_2(\chi),\,
       \int_X\varphi_+G_2(\chi_\varphi)
      \right\},
 \qquad \varphi_\pm=\max\{\pm\varphi,0\}.
\end{equation}
\end{thm}

\section{Solvability implies coerciveness}
\label{sec:solvability-coerciveness}
In this section, we prove the $(1)\implies(2)$ direction under the strict proper position condition. By Proposition 3.1 of \cite{ChenGhoshWangCoerciveness}, let
$m_i$ be the order of vanishing of $f_i$ at $r(f_2)$ for $i=1,2$. Then
\[
 0\le m_1<m_2\le d,\qquad f_1^{(m_1)}(r(f_2))<0.
\]
To simplify notations, we use the polynomial $f_i(x+r(f_2))$ and the class $\chi-r(f_2)\omega$. Then we can without loss of generality, assume that $r(f_2)=0$.

We can also extract useful consequences from the proof of Proposition 3.1 of \cite{ChenGhoshWangCoerciveness}. First, Equation~(4) of \cite{ChenGhoshWangCoerciveness} gives an $\eta>0$ such that
\begin{equation}\label{eq:strict-monomial-perturbation}
 f_2\prec f_1+\eta x^{m_1}.
\end{equation}
Second, the relation $x^{m_2}\vartriangleleft f_2$ in the same proof gives, for
$\lambda\in\Upsilon_{g_2}$,
\begin{equation}\label{eq:shifted-hessian-cone}
 \lambda \in \Gamma_{m_2}:=\Upsilon_{\sigma_{m_2}}, \qquad
 \partial_{\lambda_i}\!\left(\frac{\sigma_{m_2}}{g_2}\right)\le0.
\end{equation}

Chen-Ghosh-Wang \cite{ChenGhoshWangCoerciveness} used Aubin's $I-J$ functional $\mathcal{A}(\varphi)$ to show that $J_1(\varphi)-\varepsilon\mathcal{A}(\varphi)$ has a lower bound, and $\mathcal{A}(\varphi)\ge \delta d_1(0,\varphi)-C$ for all $\varphi\in\HH_2$ such that $I_2(\varphi)=0$. Now we want to seek for more functionals. For $0\le j\le m_2$, we first recall the Hessian-energy functional in Definition 6.1 of \cite{LuNguyenDegenerateHessian}
\[
 E_j(\varphi)=\frac{1}{j+1}\sum_{\ell=0}^j
 \int_X\varphi\,\chi_\varphi^\ell\wedge\chi^{m_2-\ell}\wedge\omega^{n-m_2}.
\]
For any smooth path $\varphi_t\in\HH_2$, integration by parts gives
\[
 \frac{d}{dt}E_j(\varphi_t)=
 \int_X\dot\varphi_t\,\chi_{\varphi_t}^j\wedge\chi^{m_2-j}\wedge\omega^{n-m_2}.
\]
We would like to use $E_{m_1}(\varphi)$ to play the role of $\mathcal{A}(\varphi)$, but it is hard to prove the coerciveness $E_{m_1}(\varphi)\ge \delta d_1(0,\varphi)-C$. Fortunately, by \eqref{eq:shifted-hessian-cone}, $E_{m_2}(\varphi)$ is convex along geodesics in $\HH_2$, and we use it to modify our functional.
Let
\[
A_{\ell}(\varphi)=\int_X\varphi\,\chi_\varphi^\ell\wedge\chi^{m_2-\ell}\wedge\omega^{n-m_2}
\]
be the $\ell$-th term in the $E_j$ functional, then
\[
A_{\ell}(\varphi)-A_{\ell+1}(\varphi)= \int_X\sqrt{-1}\partial\varphi\wedge\bar\partial\varphi\wedge
       \chi_\varphi^\ell\wedge\chi^{m_2-1-\ell}\wedge\omega^{n-m_2}\ge 0.
\]
Thus, for \[E_{m_1}=\frac{1}{m_1+1}\sum_{\ell=0}^{m_1}A_{\ell}(\varphi), \qquad E_{m_2}=\frac{1}{m_2+1}\sum_{\ell=0}^{m_2}A_{\ell}(\varphi),\]
if we want to find $\kappa \ge 0$ such that 
\[
\frac{1}{m_1+1}\sum_{\ell=1}^{m_1}A_{\ell}(\varphi)\ge \frac{\kappa}{m_2+1}\sum_{\ell=1}^{m_2}A_{\ell}(\varphi),
\]
then we can choose $\kappa=\frac{m_1(m_2+1)}{m_2(m_1+1)}<1$. It will imply that \begin{equation}\label{eq:Coercive-Em1}
E_{m_1}(\varphi)-\kappa E_{m_2}(\varphi)\ge (1-\kappa) A_0(\varphi).
\end{equation}
Then we show that $A_0$ is coercive by the following way. We start with a lemma:
\begin{lem}
\label{Laplacian-lower-bound}
For any $\lambda\in \Upsilon_{g_2}$, we have $\sum_{i=1}^{n}\lambda_i> 0$.
\end{lem}
\begin{proof}
By \eqref{eq:shifted-hessian-cone},
\[\Upsilon_{g_2}\subset \Gamma_{m_2}=\Upsilon_{\sigma_{m_2}}=\{\lambda\in\RR^n| \sigma_i >0, \forall i=1,2,...,m_2\}.\]
\end{proof}

\begin{cor}\label{lem:normalized-L1}
There is a constant $C\ge0$ such that every $u\in\HH_2$ with $\sup_Xu=0$ satisfies
\begin{equation}\label{eq:normalized-L1}
 \int_X(-u)\,\omega^n\le C, \qquad 0\le-\int_Xu\,\chi^{m_2}\wedge \omega^{n-m_2}\le C,
 \qquad 0\le-\int_Xu\,G_2(\chi)\le C.
\end{equation}
\end{cor}

\begin{proof}
By Lemma \ref{Laplacian-lower-bound},
\[
0<\operatorname{tr}_\omega\chi_u=\operatorname{tr}_\omega\chi+\Delta_\omega u.
\]
The Green-function estimate implies the upper bound on $\int_X(-u)\,\omega^n$. The other inequalities are trivial because both $\chi^{m_2}\wedge \omega^{n-m_2}$ and $G_2(\chi)$ are fixed smooth positive volume forms.
\end{proof}

\begin{prop}
There exists $\delta>0$ and $C\ge 0$ such that for any $\varphi\in \HH_2$ with $I_2(\varphi)=0$, we have
$A_0(\varphi)\ge \delta d_1(0,\varphi)-C$.
\end{prop}
\begin{proof}
If $I_2(\varphi)=0$, define $u=\varphi-\sup_X\varphi$. By \eqref{eq:normalized-L1}, $A_0(u)\ge -C$, so
\[A_0(\varphi)= (\int_X \chi^{m_2}\wedge\omega^{n-m_2}) \sup_X\varphi+A_0(u)\ge (\int_X \chi^{m_2}\wedge\omega^{n-m_2})\sup_X\varphi-C.\]
Then we apply the same argument of Lemma 9.3 of \cite{ChenGhoshWangCoerciveness} to show that $\sup_X\varphi$ is coercive in terms of $d_1$ when $I_2(\varphi)=0$. The only difference is that we use \eqref{eq:normalized-L1} here instead of the $\chi$-psh condition there.
\end{proof}

Next, as in \cite{ChenGhoshWangCoerciveness}, for $\varepsilon>0$ to be determined, we define
\begin{equation}\label{eq:perturbed-functional}
\begin{split}
 &\mathcal J_\varepsilon(\varphi)
    =J_1(\varphi)-\varepsilon (E_{m_1}(\varphi)-\kappa E_{m_2}(\varphi)), \\
 &G_{1,\varepsilon}(\chi_\varphi)
    =G_1(\chi_\varphi)-\varepsilon \bigl(\kappa\chi_\varphi^{m_2}-\chi_\varphi^{m_1}\wedge\chi^{m_2-m_1}\bigr)
       \wedge\omega^{n-m_2}, \\
 &F_\varepsilon(\chi_\varphi)=\frac{G_{1,\varepsilon}(\chi_\varphi)}{G_2(\chi_\varphi)}.
 \end{split}
\end{equation}
Define the formal derivative of $G_{1,\varepsilon}(\chi_\varphi)$ by
\[G_{1,\varepsilon}'(\chi_\varphi)=G_1'(\chi_\varphi)-\varepsilon \bigl(\kappa m_2\chi_\varphi^{m_2-1}-m_1\chi_\varphi^{m_1-1}\wedge\chi^{m_2-m_1}\bigr)
       \wedge\omega^{n-m_2},\]
where we use the convention that $0\chi_{\varphi}^{-1}=0$, even though $\chi_{\varphi}^{-1}$ is not well-defined.
Then the similar calculation in \cite{ChenGhoshWangCoerciveness} shows that
\begin{equation}
\label{eq:solution-based-first-variation}
        \frac{d}{dt}\mathcal{J}_{\varepsilon}(\varphi_t) =-\int_{X}\frac{d\varphi_t}{dt} G_{1,\varepsilon}(\chi_{\varphi_t}).
\end{equation}

Now we prove the following key lemma as an analogue of Lemma 9.4 of \cite{ChenGhoshWangCoerciveness}.
\begin{lem}\label{lem:quotient-positivity}
For all sufficiently small $\varepsilon$ and every $\varphi\in\mathcal H_2$,
\begin{equation}\label{eq:solution-perturbation-positivity}
 G_{1,\varepsilon}'(\chi_\varphi)
       -F_\varepsilon(\chi_\varphi)G_2'(\chi_\varphi) \ge 0,
 \qquad
 F_\varepsilon(\chi_\varphi)\le1.
\end{equation}
\end{lem}
\begin{proof}
The first inequality in \eqref{eq:solution-perturbation-positivity} means that its product with any test form $\theta=\sqrt{-1}\xi\wedge\bar\xi$ for a non-zero 1-form $\xi$ is non-negative.

We first deal with part of it. By \eqref{eq:strict-monomial-perturbation}, choose $\eta>0$ such that $f_2\prec f_1+\eta x^{m_1}$. By \eqref{eq:shifted-hessian-cone}, for any $\varepsilon>0$, for any $\alpha\in\Upsilon_{g_2}$,
\[
 \frac{\Pol_n(f_1+\eta x^{m_1})(\alpha)-\varepsilon \kappa \alpha^{m_2}\wedge\omega^{n-m_2}}{\Pol_n(f_2)(\alpha)}=\frac{G_1(\alpha)+\eta\,\alpha^{m_1}\wedge \omega^{n-m_1}-\varepsilon \kappa \alpha^{m_2}\wedge\omega^{n-m_2}}{G_2(\alpha)}
\]
is nondecreasing when a nonnegative $(1,1)$-form is added to $\alpha$, since the eigenvalues are nondecreasing by the Courant-Fischer-Weyl min-max principle.

Thus, it suffices to show that for all $\alpha\in\Upsilon_{g_2}(\omega)$ and all small enough $\varepsilon$, for the top form
\[G_{\eta,\varepsilon}(\alpha)=\alpha^{m_1}\wedge\bigl(\varepsilon\chi^{m_2-m_1}
       \wedge\omega^{n-m_2}-\eta\,\omega^{n-m_1}\bigr),\]
we have \[\left.\frac{d}{dt}\right|_{t=0}\frac{G_{\eta,\varepsilon}(\alpha+t\theta)}{G_2(\alpha+t\theta)}\ge 0.\]
Since $\theta^2=0$, both $G_{\eta,\varepsilon}(\alpha+t\theta)$ and $G_2(\alpha+t\theta)$ are affine in $t$, with \[
\left.\frac{d}{dt}\right|_{t=0} G_2(\alpha+t\theta)=\theta\wedge G_2'(\alpha)>0.\] In general, if $c>0$, then for $t_*=-d/c$ such that $ct_*+d=0$, we have 
\[(\frac{at+b}{ct+d})'=\frac{ad-bc}{(ct+d)^2}\ge 0\] for all $t\not=t_*$ if and only if $ad-bc=-c(at_*+b)\ge 0$. Thus, it suffices to show that for $\alpha_*=\alpha+t_*\theta\in \partial\Upsilon_{g_2}$, we have $G_{\eta,\varepsilon}(\alpha_*)\le 0$. 

To see this, let $s=(\frac{\varepsilon}{\eta})^{\frac{1}{m_2-m_1}}$. By compactness of $X$, for all $\varepsilon$ small enough, $\omega-s\chi$ is positive definite
on $X$. Then
\[
\begin{aligned}
 -G_{\eta,\varepsilon}(\alpha_*)
 &=\eta\,\alpha_*^{m_1}\wedge\bigl(\omega^{m_2-m_1}-(s\chi)^{m_2-m_1}\bigr)
       \wedge\omega^{n-m_2}\\
 &=\eta\sum_{k=1}^{m_2-m_1}\alpha_*^{m_1}\wedge(\omega-s\chi)
       \wedge\omega^{m_2-m_1-k}\wedge(s\chi)^{k-1}\wedge\omega^{n-m_2}.
\end{aligned}
\]
By \eqref{eq:shifted-hessian-cone}, $\chi\in\Gamma_{m_2}$ since $0\in\HH_2$,
and $\alpha_*\in\overline{\Gamma}_{m_2}$ by continuity. Since $\alpha_1\wedge...\alpha_{m_2}\wedge\omega^{n-m_2}\ge 0$ for all $\alpha_i\in\overline{\Gamma}_{m_2}$ by Proposition 2.1 of \cite{BlockiWeakSolution}, we see that $G_{\eta,\varepsilon}(\alpha_*)\le 0$.

Finally, $F_\varepsilon(\alpha)\le \lim_{t\to\infty}F_\varepsilon(\alpha+t\omega)=1-\varepsilon\kappa\le 1$ if $m_2=d$, and the limit is $1$ if $m_2<d$.
\end{proof}

\begin{proof}[Proof of $(1) \implies (2)$ in Theorem~\ref{thm:main}]
Assume that $\chi_{\varphi_*}$ satisfies $G_1(\chi_{\varphi_*})=0$. Then we solve the equation $G_{1,\varepsilon}(\chi_\varphi)=C_* \varepsilon G_2(\chi_\varphi)$ for the constant $C_*$ such that \[
\int_X G_{1,\varepsilon}(\chi)=C_* \varepsilon \int_X G_2(\chi)\]
for all $\varepsilon>0$ small. We view it as a perturbation of the equation $G_1(\chi_{\varphi_*})=0$. The linearization in $\varphi$ is elliptic, self-adjoint, with kernel the constants. So by the implicit function theorem, there is a solution $\varphi_{\varepsilon}$ close to $\varphi_*$. For $\varepsilon$ small, we can assume that $\varphi_{\varepsilon} \in \HH_2$. By adding a constant, we can also assume that $I_2(\varphi_{\varepsilon})=0$.

Fix one such $\varepsilon>0$. For each $\varphi\in\HH_2$ with $I_2(\varphi)=0$, we take the
$\rho$-geodesic $\varphi_{\rho,t}$ from $\varphi_\varepsilon$ to $\varphi$ in Theorem~\ref{prop:geodesics}. Define $\mathcal J_\rho(t)=\mathcal J_\varepsilon(\varphi_{\rho,t})+C_*\varepsilon I_2(\varphi_{\rho,t})$. Then $\mathcal J'_\rho(0)=0$, and
\[
\begin{split}
 \mathcal J''_\rho(t)
 &=\int_X\sqrt{-1}\partial \dot\varphi_{\rho,t}\wedge\bar\partial \dot\varphi_{\rho,t}\wedge
       (G'_{1,\varepsilon}-F_\varepsilon G'_2)(\chi_{\varphi_{\rho,t}})
    -\int_X (F_\varepsilon-C_*\varepsilon)\,\mathcal D_\rho(t)\\
    &\ge-(1+|C_*|\varepsilon)\int_X\mathcal D_\rho(t)\ge -C\rho^{1/d},
\end{split}
\]
where
$0\le \mathcal D_\rho\le C\rho^{1/d}\omega^n$ is the defect term defined in \eqref{eq:geodesic-defect}.
Thus,
\[
 \mathcal J_\varepsilon(\varphi)-\mathcal J_\varepsilon(\varphi_{\varepsilon}) = \int_0^1(1-t)\mathcal J_\rho''(t)\,dt\ge -C\rho^{1/d}.
\]
Letting $\rho\to0$, we obtain
$\mathcal J_\varepsilon(\varphi)\ge\mathcal J_\varepsilon(\varphi_{\varepsilon})$. This finishes the proof.
\end{proof}

\section{Coerciveness implies numerical criterion}
\label{sec:numerical}
In this section, we prove the $(2)\implies (3)$ direction in Theorem \ref{thm:main}. By Lemma 2.22 of \cite{ChenNieXuNumericalCriterion}, this holds for $G_2$. In other words, using the assumption that $0\in\HH_2$, we have
\begin{equation}\label{eq:positive-denominator-intersection}
 \int_Y G_2^{(r)}(\chi)>0
\end{equation}
for every irreducible subvariety $Y$ with codimension $r\le d$. We need to prove the similar thing for $G_1$.

\subsection{Energy slopes of logarithmic regularizations}

\begin{lem}\label{lem:logarithmic-slopes}
Fix an irreducible reduced analytic subvariety $Y$ of codimension
$1\le r\le d$. There are smooth functions $v_T\le0$, $T\ge1$,
with $\ddc v_T\ge-A\omega$ and $\int_X(-v_T)\omega^n\le A$,
where $A$ is independent of $T$, such that the following holds.
For all sufficiently small $a>0$, $a v_T$ belongs to
$\HH_2$. For any $i=1,2$,
\begin{equation}\label{eq:logarithmic-energy-slope}
 \lim_{T\to\infty}\frac{-I_i(a v_T)}{T}
   =\frac{(2\pi)^r}{(r+1)!} a^{r+1}\int_Y G_i^{(r)}(\chi)+O(a^{r+2}).
\end{equation}
\end{lem}

\begin{proof}
Choose a finite coordinate cover $\{U_\alpha\}$ such that the reduced ideal $\mathcal I_Y$ has holomorphic generators $f_{\alpha1},\ldots,f_{\alpha N_\alpha}$ on a neighborhood of each $\overline{U_\alpha}$. Choose real smooth functions $\eta_\alpha$ with compact support in $U_\alpha$ and $\sum_\alpha\eta_\alpha^2=1$ on $X$. Set
\begin{equation}\label{eq:global-log-model}
 v_T=\log(\sum_\alpha\eta_\alpha^2\sum_j|f_{\alpha j}|^2+e^{-T})-C_0,
\end{equation}
where $C_0\ge\max\{1,\log(\sup_X\sum_\alpha\eta_\alpha^2\sum_j|f_{\alpha j}|^2+e^{-1})\}$.
Then $-T-C_0\le v_T\le0$, and $v_\infty=\lim_{T\to\infty} v_T$ satisfies $v_\infty^{-1}(-\infty)=Y$.

By the calculation of Demailly-P\u{a}un \cite[Lemma~2.1(ii), equations~(2.2)-(2.3), p.~1255]{DemaillyPaunKahlerCone}, there exists $A$ independent of $T$ such that
\[\ddc v_T\ge-A\omega.\]
Moreover, $\displaystyle{\log \sum_\alpha\eta_\alpha^2\sum_j|f_{\alpha j}|^2\in L^1(X)}$ and \[0\le-v_T\le C_0-\log \sum_\alpha\eta_\alpha^2\sum_j|f_{\alpha j}|^2\quad\text{a.e. on }X.\]
Thus, by enlarging $A$, we can assume that $\int_X(-v_T)\omega^n\le A$ also holds.

Since $X$ is compact and $\lambda(\omega^{-1}\chi)\in\Upsilon_{g_2}$, there exists $a_0>0$ such that
\begin{equation}\label{eq:small-shift-admissible}
 \lambda\bigl(\omega^{-1}(\chi-aA\omega)\bigr)\in\Upsilon_{g_2}
 \qquad(0\le a\le a_0).
\end{equation}
From the decomposition
\[
 \chi_{a v_T}=(\chi-aA\omega)+a(\ddc v_T+A\omega),\qquad
 a(\ddc v_T+A\omega)\ge0,
\]
we get $a v_T\in\HH_2$.

Take a principalization $\pi:\widetilde X\to X$ of $\mathcal I_Y$:
\[
 \mathcal I_Y\mathcal O_{\widetilde X}=\mathcal O_{\widetilde X}(-E),
 \qquad E=\sum_\nu m_\nu E_\nu.
\] For the calculation below, we take the following one: first resolve the singularity of $Y$ by $\pi_1: X_1\to X$, such that the strict transform $Y_1$ of $Y$ is smooth in $X_1$, $D\subset \pi_1^{-1}(Y_{\rm sing})$ and  $\mathcal{I}_Y\mathcal O_{X_1}=\mathcal O_{X_1}(-D)\mathcal I_{Y_1}$. $\pi_1$ is biholomorphic on $\pi_1^{-1}(X\setminus Y_{\rm sing})$. Then we blow up along $Y_1$, and get $\pi_2:X_2\to X_1$. Let $F=\pi_2^{-1}(Y_1)$. Then $\mathcal I_Y\mathcal O_{X_2}=\mathcal O_{X_2}(-\pi_2^{*}D)\mathcal O_{X_2}(-F)=\mathcal O_{X_2}(-F-\pi_2^{*}D).$  Take $\widetilde{X}=X_2$ and $E=F+\pi_2^{*}D$. 

Since the pulled-back generators generate the principal ideal of $E$, there is a smooth Hermitian metric $h_E$ on $\mathcal O(E)$ such that
\[
 \sum_\alpha(\eta_\alpha\circ \pi)^2\sum_j|f_{\alpha j}\circ\pi|^2=|s_E|_{h_E}^2,\qquad \ddc\log|s_E|_{h_E}^2=-\theta_E+2\pi[E],
\]
where $\theta_E$ is the Chern curvature of $(\mathcal{O}(E),h_E)$.
Set
\[
 \widetilde\chi=\pi^*\chi,\quad \widetilde\omega=\pi^*\omega,\quad
 \widetilde v_T=v_T\circ\pi=\log(|s_E|_{h_E}^2+e^{-T})-C_0.
\]
The calculations in \cite[\S5, proof of Lemma~9, Step~3, equations~(16)-(21)]{LiuSemistableJEquation} give, for any $s\geq 0$ and $0\le p\le d$, \begin{equation}\label{eq:liu-monomial-slope}
\begin{aligned}
 &\lim_{T\to\infty}\int_{\widetilde X}\frac{-\widetilde v_T}{T}(\widetilde\chi+s\ddc\widetilde v_T)^p\wedge\widetilde\omega^{n-p}=\int_{\widetilde X}\bigl[\widetilde\chi^p-(\widetilde\chi-s\theta_E)^p\bigr]\wedge\widetilde\omega^{n-p}.
\end{aligned}
\end{equation}

By \eqref{eq:energy-variation},
\begin{equation}
 I_i(a v_T)
     =\int_0^1 \Big(\int_X av_T\,G_i(\chi_{tav_T})\Big)dt.
\end{equation}
Thus
\begin{equation}\label{eq:exact-slope-polynomial}
\begin{aligned}
 \lim_{T\to\infty}\frac{-I_i(av_T)}{T}
 &=a\int_0^1 \Big(\int_{\widetilde X} G_i(\widetilde\chi)- G_i(\widetilde\chi-ta\theta_E)\Big)dt\\
 &=-\sum_{m=1}^d\frac{a^{m+1}}{(m+1)!}
       \int_{\widetilde X}(-\theta_E)^m\wedge\pi^*G_i^{(m)}(\chi).
\end{aligned}
\end{equation}

Next, since $[\theta_E]=2\pi [E]$, for any smooth closed $(n-m,n-m)$-form $\eta$,
\[
 \int_{\widetilde X}\theta_E^m\wedge\pi^*\eta
 =2\pi\int_E\theta_E^{m-1}\wedge\pi^*\eta.
\]
For $m<r$, the right-hand side vanishes since
$\operatorname{rank}_{\mathbb C}d(\pi|_{E_\nu})\le n-r<n-m$ generically.
For $m=r$, only components dominating $Y$ contribute. For $y\in Y_{\rm reg}$, $\pi^{-1}(y)=\mathbb{P}^{r-1}$, and 
\[
 \mathcal O(E)|_{\mathbb P^{r-1}}=\mathcal O_{\mathbb P^{r-1}}(-1),
 \qquad \int_{\mathbb P^{r-1}}(-\theta_E)^{r-1}=(2\pi)^{r-1}.
\]
Hence \[
\int_{E}\theta_E^{r-1}\wedge\pi^*\eta=\int_{Y_{\rm reg}}\left(\int_{\pi^{-1}(y)}\theta_E^{r-1}\right)\eta(y)=(-2\pi)^{r-1}\int_{Y_{\rm reg}}\eta.
\]
For $r=1$, this is the integral of $1$ over a point. Thus
\begin{equation}\label{eq:leading-coefficients}
 \int_{\widetilde X}(-\theta_E)^m\wedge\pi^*G_i^{(m)}(\chi)
 =\begin{cases}
 0,&1\le m<r,\\
 -(2\pi)^r \int_Y G_i^{(r)}(\chi),&m=r.
 \end{cases}
\end{equation}
Substitution into \eqref{eq:exact-slope-polynomial} yields
\begin{equation}\label{eq:slope-remainder}
\begin{aligned}
 \lim_{T\to\infty}\frac{-I_i(av_T)}{T}&=\frac{(2\pi)^r}{(r+1)!}a^{r+1}\int_Y G_i^{(r)}(\chi)-\sum_{m=r+1}^d\frac{a^{m+1}}{(m+1)!}
       \int_{\widetilde X}(-\theta_E)^m\wedge\pi^*G_i^{(m)}(\chi).
\end{aligned}
\end{equation}
After decreasing $a_0$ so that $a_0\le1$, we use the fixed coefficients of this finite polynomial to obtain
\[
 |\sum_{m=r+1}^d\frac{a^{m+1}}{(m+1)!}
       \int_{\widetilde X}(-\theta_E)^m\wedge\pi^*G_i^{(m)}(\chi)|\le \sum_{m=r+1}^d\frac{a^{r+2}}{(m+1)!}
       \left|\int_{\widetilde X}(-\theta_E)^m\wedge\pi^*G_i^{(m)}(\chi)\right|
\]
for all $0<a\le a_0$.
This now proves \eqref{eq:logarithmic-energy-slope}.

\end{proof}

\subsection{The distance before and after normalization}

\begin{lem}\label{lem:logarithmic-distance}
Set 
\begin{equation}\label{eq:normalized-logarithmic-test}
\varphi_{a,T}=a v_T-\frac{I_2(a v_T)}{V_2}.
\end{equation}
Then $I_2(\varphi_{a,T})=0$. Moreover, there are constants
$K,C_a\ge0$, independent of $T$, such that, for small $a$,
\begin{equation}\label{eq:normalized-log-distance}
 d_1(0,\varphi_{a,T}) \ge -(1-Ka)I_2(a v_T)-C_a.
\end{equation}
\end{lem}

\begin{proof}
Choose $a_0$ as in \eqref{eq:small-shift-admissible}. Set
\begin{equation}\label{eq:comparison-constant}
 K=\sup_{\substack{x\in X\\0\le s\le a_0}}
 \frac{A\omega\wedge G'_2(\chi-sA\omega)}{G_2(\chi)}<\infty.
\end{equation}
Then we have
\begin{equation}\label{eq:shifted-volume}
 G_2(\chi-aA\omega)=G_2(\chi)-\int_0^a A\omega\wedge G'_2(\chi-sA\omega)\,ds\ge(1-Ka)G_2(\chi).
\end{equation}
Decrease $a_0$ if necessary so that $Ka\le1/2$. Since
\[\chi_{a v_T}-(\chi-aA\omega)=a(\ddc v_T+A\omega)\ge0,\]
the monotonicity of $G_2$ gives
\begin{equation}\label{eq:test-volume-comparison}
 G_2(\chi_{\varphi_{a,T}})=G_2(\chi_{a v_T})
 \ge G_2(\chi-aA\omega)\ge(1-Ka)G_2(\chi).
\end{equation}
Now apply the second term in Theorem~\ref{ass:distance}:
\begin{align}
 d_1(0,\varphi_{a,T})
 &\ge\int_X(\varphi_{a,T})_+G_2(\chi_{\varphi_{a,T}})\ge(1-Ka)\int_X\varphi_{a,T}G_2(\chi)\notag\\
 &=(1-Ka)\left(-I_2(a v_T)+\int_X a v_T G_2(\chi)\right).\label{eq:normalized-distance-lower}
\end{align}
If $M_2=\sup_XG_2(\chi)/\omega^n$, then we have 
\[0\le-\int_X a v_T G_2(\chi)
 \le aM_2\int_X(-v_T)\omega^n\le aM_2A.\]
Thus one may take $C_a=(1-Ka)aM_2A$.
\end{proof}

\begin{prop}
\label{prop:coercivity-numerical}
Assume the coerciveness condition \eqref{eq:main-d1-coercivity}. Then every
irreducible analytic $Y\subsetneq X$ with
$1\le\operatorname{codim}Y\le d$ satisfies
\begin{equation}\label{eq:uniform-unmixed-numerical}
 \int_YG_1^{(r)}(\chi)
 \ge\delta\int_YG_2^{(r)}(\chi)>0,
 \qquad r=\operatorname{codim}Y.
\end{equation}
Here $\delta$ is the coerciveness constant for the distance
normalization \eqref{eq:d1-definition}.
\end{prop}

\begin{proof}Fix $Y$ and use its tests from Lemma~\ref{lem:logarithmic-slopes}. 
Since $V_1=0$, \[
J_1(\varphi_{a,T})=-I_1(a v_T-\frac{I_2(a v_T)}{V_2})=-I_1(a v_T)=J_1(a v_T).\]
Apply coerciveness to \eqref{eq:normalized-logarithmic-test} and
then use Lemma~\ref{lem:logarithmic-distance}:
\[
 J_1(a v_T)\ge\delta(1-Ka)(-I_2(a v_T))-\delta C_a-C.
\]
Divide by $T$ and apply Lemma~\ref{lem:logarithmic-slopes}. First let $T\to\infty$, then divide by
$\frac{(2\pi)^r}{(r+1)!}a^{r+1}$ and let $a\downarrow0$. This gives
\eqref{eq:uniform-unmixed-numerical}.
\end{proof}

\section{Numerical criterion implies solvability}
\label{sec:continuity}

In this section, we prove the $(3)\implies (1)$ direction in Theorem \ref{thm:main}. Keep $(X,\omega,\chi,f_2)$ fixed and equip the space of real
polynomials of degree at most $d$ with its coefficient topology.
Let $\mathcal P_2$ be the set of monic degree-$d$ right-Noetherian
polynomials $f$ such that
\[
 f_2\prec f\ \text{strictly},\qquad
 \int_XG(\chi)=0,
\]
where $g$ and $G$ are the polarizations of $f$.

We first show that $\mathcal P_2$ is convex. To see this, if $f_{1,0}, f_{1,1} \in \mathcal P_2$, then for 
\[
f_{1,t}=tf_{1,1}+(1-t)f_{1,0},
\] by Theorem 6.15 (vi) of \cite{FangMaGarding}, $f_{1,t}$ is also right-Noetherian and we have $f_2\prec f_{1,t}$ for all $0\le t\le1$. Using \[\partial_{\lambda_i}\frac{g_{1,t}}{g_2}=t\partial_{\lambda_i}\frac{g_{1,1}}{g_2}+(1-t)\partial_{\lambda_i}\frac{g_{1,0}}{g_2}>0\] and Proposition 3.1 of \cite{ChenGhoshWangCoerciveness}, we have $f_2\prec f_{1,t}$ strictly for all $0\le t\le1$. So $f_{1,t}\in \mathcal P_2$ since \[
\int_XG_{1,t}(\chi)=t\int_XG_{1,1}(\chi)+(1-t)\int_XG_{1,0}(\chi)=0.
\]

Define the convex subset
\[
 \mathcal N_2
 =\{f\in\mathcal P_2:\int_Y G^{(r)}(\chi)>0
       \text{ for every irreducible }Y\subsetneq X,
       \ 1\le r=\operatorname{codim}Y\le d\},
\]
and the subset
\[
\mathcal S_2
  =\{f\in\mathcal P_2:\text{there exists }\chi_\varphi \in \partial\Upsilon_{G}\subset\Upsilon_{G_2}\}\subset \mathcal N_2.
\]

\begin{thm}
\label{thm:connected-component}
The solvable locus $\mathcal S_2=\emptyset$ or $\mathcal S_2=\mathcal N_2$.
\end{thm}

Using
the numerical criteria of Fu-Li-Zhang-Zhang
\cite{FuLiZhangZhangNumerical} and the spectral solvability theorem
\cite[Theorem~1.1]{FangMaSpectral}, it suffices to prove the positivity of the mixed integral. As a useful lemma, the log concavity inequality below comes from its pointwise version in Fang and Ma \cite{FangMaIdealGarding}, and was used in Section 10 of \cite{FuLiZhangZhangNumerical}.

\begin{lem}
\label{LogConcaveLemma}

Let $f$ be a monic right-Noetherian polynomial with degree $2\le d\le n$, and define $g,G$ as in \eqref{eq:polynomials-forms}. Suppose that there exists $\chi_\varphi$ such that $\lambda(\omega^{-1}\chi_\varphi)\in \Upsilon_{g}$ at each point. Let $Y$ be irreducible of dimension $p$. For $0\leq r \leq p-n+d$, we define 
\[
U_r = \int_Y \omega^{p-n+d-r}\wedge G^{(d-r)}(\chi)
\]
Then $U_r > 0$, and for any $2 \leq r \leq p-n+d$, we have $U_{r-1}^2 \geq U_{r-2}U_{r} $.
\end{lem}

\begin{proof}
    Since this is purely homological, without loss of generality, we assume $\varphi=0$. 
    Then $G^{(d-r)}(\chi)$ is a positive form on $X$. So we have $U_r>0$.
    
    For the second statement, first we assume that $Y$ is smooth.  
    Fix $2\leq r\leq p-n+d$, and take the constant $c= \frac{U_{r-1}}{U_{r}}$.
    Assume that $f(x) = \sum_{k=0}^{d}a_k\binom{n}{k}x^k$, where $a_d = \binom{n}{d}^{-1}.$ Then 
    \[
    G^{(d-r)}(\chi) = \frac{n!}{(n-d+r)!}\sum_{k=0}^{r}a_{k+d-r}\binom{n-d+r}{k}\chi^{k}\wedge\omega^{n-d+r-k}.
    \]
    Then we have
    \begin{align*}
    \omega^{p-n+d-r}\wedge G^{(d-r)}(\chi)=\frac{n!}{(n-d+r)!}\sum_{k=0}^{r}a_{k+d-r}\binom{n-d+r}{k}\chi^{k}\wedge\omega^{p-k}.
    \end{align*}
    For simplicity, we write $\tilde{\lambda}(\omega^{-1}\chi)\in\mathbb{R}^p$ as the eigenvalues of $(\omega|_{Y})^{-1}(\chi|_Y)$, and
    \[
    q(\tilde{\lambda}) = \Pol_{p}(f^{(d-r)})(\tilde{\lambda}) = \frac{n!}{(n-d+r)!}\sum_{k=0}^{r}a_{k+d-r}\binom{n-d+r}{k}\frac{\sigma_{k}(\tilde{\lambda})}{\binom{p}{k}},
    \]
    and we denote $Q(\chi) = q(\tilde{\lambda}(\omega^{-1}\chi)) = \sum_{k=0}^{r}b_k\sigma_k(\tilde\lambda),$ where $b_k =a_{k+d-r}\frac{n!}{(n-d+r)!}\binom{n-d+r}{k}\binom{p}{k}^{-1}. $
    Then $U_r = \int_Y Q(\chi)\omega^p$. Moreover, by Lemma 2.22 of \cite{ChenNieXuNumericalCriterion}, $G^{(i)}(\chi)$ is positive for any $n-p\le i \le d$. So $G^{(i)}(\chi)|_{Y}$ is also positive. Then for $2\leq r \leq p-n+d$, and any $t>0$, \[Q(\chi+t\omega)\omega^p|_{Y} = \left. G^{(d-r)}(\chi+t\omega)\wedge\omega^{d-r-n+p} \right|_Y
    =\sum_{j=0}^{r}\frac{t^j}{j!} \left. G^{(d-r+j)}(\chi)\wedge\omega^{d-r-n+p+j} \right|_Y
    \] is also positive, this implies that $\tilde{\lambda}(\omega^{-1}\chi) \in \Upsilon_Q$.

    Define 
    \[
    DQ_{\chi}(\alpha)=\frac{d}{dt}\Big|_{t=0}Q(\chi+t\alpha), \qquad D^2Q_{\chi}(\alpha)=\frac{d^2}{dt^2}\Big|_{t=0}Q(\chi+t\alpha).
    \]
    By computation,
    \[
     DQ_{\chi}(\alpha)\omega^p=\alpha\wedge\sum_{k=1}^{r}kb_k\binom{p}{k}\chi^{k-1}\wedge\omega^{p-k},
    \]
    and we notice that
    \[
    \int_Y DQ_{\chi}(\omega + \ddc u)\omega^p = \int_Y\omega\wedge \sum_{k=1}^{r}kb_k\binom{p}{k}\chi^{k-1}\wedge\omega^{p-k} = U_{r-1}.
    \]

    Since $f$ is right-Noetherian and its polarization is ideal Fang-Ma-G{\aa}rding, the same holds for $f^{(d-r)}$ and $q(\tilde{\lambda}) = \Pol_{p}(f^{(d-r)}).$
    Using Lemma 2.5 and Lemma 2.7 of Fang-Ma\cite{FangMaSpectral}, we obtain that $Q$ is an ideal G{\aa}rding elliptic operator and $\log Q$ is concave on $\Upsilon_Q$. In particular, for any smooth $u$,
    \[
    D^2Q_{\chi}(\omega+\ddc u)Q(\chi) \leq DQ_{\chi} (\omega+\ddc u)^2.
    \]
    Since $DQ_\chi$ is elliptic and self-adjoint, and its kernel consists of constants, there exists a smooth function $u$ such that
    \begin{equation}
    \label{LogConcave-Equation}
    DQ_{\chi} (\omega+\ddc u) = c Q(\chi),
    \end{equation}
    where $c= \frac{U_{r-1}}{U_{r}}$ is the constant such that the integral of \eqref{LogConcave-Equation} on $Y$ matches.
    Then
    \begin{align*}
        U_{r-2} &= \int_Y  D^2Q_{\chi}(\omega+\ddc u)\omega^p\leq \int_Y\frac{ DQ_{\chi}(\omega+\ddc u)^2}{Q(\chi)}\omega^p\\
                &= \int_Y cDQ_{\chi}(\omega+\ddc u)\omega^p =c U_{r-1},
    \end{align*}
    which proves that $U_{r-2}U_{r}\leq U_{r-1}^2.$

    For singular $Y$, choose a compact K{\"a}hler resolution $\pi:\widetilde X\to X$ such that the proper transform $\tilde Y$ of $Y$ is smooth. Fix a K{\"a}hler form $\eta$ on $\widetilde X$. For $\varepsilon>0$, use
    \[
     \tilde{\omega}_\varepsilon=\pi^*\omega+\varepsilon\eta,\qquad
     \tilde{\chi}_\varepsilon=\pi^*\chi+C\varepsilon\eta,
    \]
    where $C$ is a constant independent of $\varepsilon$ such that $C\omega>\chi$. On the open set where \(\pi\) is an isomorphism onto the smooth locus,  we have
    \[
    C\tilde{\omega}_{\varepsilon} - \tilde{\chi}_{\varepsilon} = C\pi^*\omega-\pi^*\chi \ge 0
    \]
    and
    \[
    \tilde{\omega}_{\varepsilon}^{-1}\tilde{\chi}_{\varepsilon} = CI - \tilde{\omega}_{\varepsilon}^{-1}(C\pi^*\omega-\pi^*\chi).
    \]
    So for any $\varepsilon>0$, by the min-max principle, the $i$-th largest eigenvalue satisfies
    \[\lambda_i(\tilde{\omega}_{\varepsilon}^{-1}\tilde{\chi}_{\varepsilon})=C-\min_{V\subset\mathbb{C}^n, \dim V=i}\max_{0\neq v\in V}\frac{(C\pi^*\omega-\pi^*\chi)(v,\bar v)}{\tilde{\omega}_{\varepsilon}(v,\bar v)}.
    \]
    Therefore, we have
    \[ \lambda(\tilde{\omega}_{\varepsilon}^{-1}\tilde{\chi}_{\varepsilon})(x) \in \lambda(\omega^{-1}\chi)(\pi(x)) + \mathbb{R}^n_{\geq 0}\subset \Upsilon_{g}
    \]
    for any $x\in \tilde X$ outside the exceptional locus, and this also holds across the exceptional locus by the continuity. Then we can apply the smooth proof to $(\widetilde X, \widetilde Y,\tilde{\omega}_\varepsilon,\tilde{\chi}_\varepsilon)$ and obtain that
    \begin{equation}
    \label{SingularLog-Concave}
    \tilde{U}_{r-2,\varepsilon}\tilde{U}_{r,\varepsilon}\leq \tilde{U}_{r-1,\varepsilon}^2.        
    \end{equation}
    Here we define
    \[
    \tilde{U}_{r,\varepsilon} = \int_{\widetilde Y} \tilde{\omega}_{\varepsilon}^{p-n+d-r}\wedge G^{(d-r)}_{\tilde{\omega}_{\varepsilon}}(\tilde{\chi}_{\varepsilon}),
    \]
    where $G^{(d-r)}_{\tilde{\omega}_{\varepsilon}}$ means that we use $\tilde{\omega}_{\varepsilon}$ to replace the form $\omega$ in $G^{(d-r)}$.

    All of the integrals $\tilde{U}_{r,\varepsilon}$ are polynomials in $\varepsilon$, so \eqref{SingularLog-Concave} converges as $\varepsilon\to 0$. Integration on the resolution agrees with integration on $Y$ for pulled-back smooth top forms. So we obtain that $U_{r-2}U_{r}\leq U_{r-1}^2.$

\end{proof}

Back to the proof of Theorem \ref{thm:connected-component}:
\begin{proof}
    Since $\mathcal N_2$ is connected, it suffices to show that $\mathcal S_2$ is both open and closed. The openness of $\mathcal S_2$ follows from the openness of the $C$-subsolution and Theorem 1.1 of \cite{FangMaSpectral}. 
    
    To prove the closedness, suppose $f_i\in\mathcal{S}_2$  and $f_i \to f\in \mathcal{N}_2.$
    By Theorem 1.2 of  Fu-Li-Zhang-Zhang \cite{FuLiZhangZhangNumerical} and Corollary 1.2 of Fang-Ma \cite{FangMaSpectral}, we only need to prove the mixed intersection inequalities with $\kappa = \omega$.
    Fixing an irreducible analytic $Y$ of dimension $p$, for $0\leq r\leq p-n+d$, we define $ U_r(f_i)$ as in Lemma \ref{LogConcaveLemma}. Then we only need to prove $U_r(f)>0$ for all $r<d$. By Lemma \ref{lem:strict-RN-and-inclusion} and its proof, $f$ is strictly right-Noetherian, and there exists $\varepsilon>0$, such that $f_i+\varepsilon$ is also right-Noetherian for large $i$. Since $\bar\Upsilon_{g_i} \subset \Upsilon_{g_i+\varepsilon}$, we can apply Lemma \ref{LogConcaveLemma} to $f_i+\varepsilon$. Let $s= p-n+d$.
    For $0\leq r < s$, Lemma \ref{LogConcaveLemma} implies that
    \[
    U_r(f_i+\varepsilon) \geq U_0^{1-r/s}(f_i+\varepsilon) U^{r/s}_{s}(f_i+\varepsilon).
    \]
    By taking the limit of $f_i$, we obtain that 
    \[
    U_r(f+\varepsilon) \geq U_0^{1-r/s}(f+\varepsilon) U^{r/s}_{s}(f+\varepsilon).
    \]
    
    Note that $U_r(f+\varepsilon)=U_r(f)$ for $r<d$, and when $r=d$, \[U_d(f+\varepsilon) = \int_X \bigl( G(\chi)+\varepsilon\omega^n \bigr)= \varepsilon\int_X\omega^n>0.\]
    Thus, together with the condition that $f\in \mathcal{N}_2$, we see that $U_s(f+\varepsilon)>0$. This finishes the proof since 
    \[U_0(f+\varepsilon) = d!\int_Y \omega^p>0.\]
    
\end{proof}

Now we show that $\mathcal S_2$ is non-empty.

\begin{lem}\label{lem:solvable-anchor}
The polynomial $f_2-v_2 \in \mathcal{S}_2$.
\end{lem}

\begin{proof}
Note that all the higher derivatives of $f_2-v_2$ are the same as the higher derivatives of $f_2$, and $f_2(r(f_2))=0$, and $f_2$ is increasing on $[r(f_2),\infty)$. Thus, $r(f_2-v_2)>r(f_2)$, and $f_2-v_2\le 0$ on $[r(f_2),r(f_2-v_2)]$. So $f_2-v_2$ is strictly right-Noetherian and $f_2\prec f_2-v_2$. Since $(f_2-v_2)(r(f_2))<0$, by Proposition 3.1 of \cite{ChenGhoshWangCoerciveness}, $f_2\prec f_2-v_2$ strictly.

Note that the polarization of $f_2-v_2$ is $G_2(\chi_\varphi)-v_2\omega^n$. The $C$-subsolution condition is the same as that for $f_2$, and the assumption $\HH_2\not=\emptyset$ shows that there exists a smooth $C$-subsolution. Together with the integral condition
$\int_X G_2(\chi)=V_2=v_2\int_X\omega^n$, Theorem 1.1 of \cite{FangMaSpectral} shows that
$f_2-v_2\in \mathcal S_2$.
\end{proof}

Thus, $f_1 \in \mathcal{N}_2=\mathcal{S}_2$ and this finishes the proof of $(3)\implies (1)$ in Theorem \ref{thm:main}.

\begingroup
\small
\bibliographystyle{plain}
\bibliography{references}

\begin{center}
\begin{minipage}{0.88\textwidth}
\small
\textbf{Gao Chen}\\[0.35em]
University of Science and Technology of China\\
No.\ 96 Jinzhai Road, Baohe District, Hefei, Anhui 230000, P.\,R.\ China\\[0.35em]
\href{mailto:chengao1@ustc.edu.cn}{\texttt{chengao1@ustc.edu.cn}}\\
\textbf{Junbang Liu}\\[0.35em]
Department of Mathematics\\
The Hong Kong University of Science and Technology\\
Clear Water Bay, Kowloon, Hong Kong\\[0.35em]
\href{mailto:junbangliu@ust.hk}{\texttt{junbangliu@ust.hk}}\\
\textbf{Sijie Nie}\\[0.35em]
University of Science and Technology of China\\
No.\ 96 Jinzhai Road, Baohe District, Hefei, Anhui 230000, P.\,R.\ China\\[0.35em]
\href{mailto:niesijie@mail.ustc.edu.cn}{\texttt{niesijie@mail.ustc.edu.cn}}
\end{minipage}
\end{center}

\endgroup

\end{document}